\documentclass{article}
 \usepackage[utf8]{inputenc}
\usepackage[T1]{fontenc}
\usepackage[english]{babel}
\usepackage{times}
\usepackage{amsmath}
\usepackage{amsthm}
\usepackage{amssymb}
\usepackage{mathrsfs}
\usepackage{graphicx}
\usepackage{setspace}
\usepackage{tikz}
\usepackage{tikz-cd}
\usepackage{nomencl}

\makeatletter 
\renewcommand\@biblabel[1]{}
\makeatother

\makenomenclature

\title{Stability of symmetric cube gamma factors for $\operatorname{GL}(2)$}
\author{Daniel Shankman and Dongming She}
\date{}

\DeclareMathOperator{\GL}{GL}

\DeclareMathOperator{\Det}{Det}

\begin{document}

\maketitle

\section*{Introduction}

Let $\mathbf M$ be a connected, reductive group over a $p$-adic field $k$.  Let $\psi$ be a nontrivial character of $k$, and let $r$ be a finite dimensional complex representation of the Langlands dual group $^L\mathbf M$ of $\mathbf M$ whose restriction to the connected component $^L\mathbf M^{\circ}$ is complex analytic.  For certain $r$, there are a handful of methods to attach local gamma factors $\gamma(s,\pi,r,\psi)$ to irreducible, admissible representations $\pi$ of $\mathbf M(k)$, for example the Langlands-Shahidi method [Sh90] or various methods of integral representations ([PiRa86], [JaPiSh83]) .  The definition of $\gamma(s,\pi,r,\psi)$ for general $\pi$ and $r$ remains conjectural.  

When the gamma factor is defined, we expect several properties to hold.  Among them is the property of \emph{stability}: that is, we expect $\gamma(s,\pi,r,\psi)$ to only depend on the central character of $\pi$, up to a highly ramified twist.  After all, once a local Langlands correspondence for $\mathbf M$ is established, the gamma factor $\gamma(s,\pi,r,\psi)$ should be equal to a corresponding Artin gamma factor, and an analogous stability property for Artin factors is known [De72].

In fact, for $\mathbf M = \GL_n$, where we do have a local Langlands correspondence ([He00], [HaTa01], [Sc13]), there is a general method using global arguments to prove that $\gamma(s,\pi,r,\psi)$ equals the corresponding Artin gamma factor under the local Langlands correspondence (Theorem 1.2 of [Sh12]).  This method was carried out successfully in the case of symmetric and exterior square representations [CoShTs17], twisted symmetric square representations [She19], and the Asai representation [Shan18].  But in order for this method to work, there are a number of properties which $\gamma(s,\pi,r,\psi)$ must be already known to have, stability among them.

One of the most successful methods of showing stability of gamma factors, at least for generic supercuspidal representations, is the method of asymptotic expansion of partial Bessel functions.  Very loosely speaking, one expresses $\gamma(s,\pi,r,\psi)$ as a Mellin transform $\operatorname{MT}$ of a partial Bessel function.  The problem then becomes to write the partial Bessel function as a sum of two functions $F_{\omega_{\pi}}$ and $F_{\operatorname{smooth}}$, the first  depending only on the central character of $\pi$, and the second having a certain ``smoothness'' property, such that if we write
\[ \gamma(s,\pi,r,\psi) = \operatorname{MT}(F_{\omega_{\pi}}) + \operatorname{MT}(F_{\operatorname{smooth}}) \]
then $\operatorname{MT}(F_{\operatorname{smooth}})$ will be zero provided the central character of $\pi$ is sufficiently highly ramified, leaving us only with a dependence on the central character of $\pi$.

This method of showing stability for supercuspidals has been successfully carried out for various gamma factors ([CoPi98], [AsSh06], [CoShTs17], [Shan18], [She19]).  In the case of Langlands-Shahidi gamma factors for maximal self-associate parabolics, there is a general method to express gamma factors as Mellin transforms of partial Bessel functions [Sh02].  But the asymptotic expansion of these partial Bessel functions is still unsolved in the general case.  

Our main result in this paper is the stability of the symmetric cube gamma factor for $\operatorname{GL}_2$, for supercuspidal representations.  This factor is defined by the Langlands-Shahidi method by embedding $\operatorname{GL}_2$ as a maximal self-associate Levi subgroup of the exceptional Lie group $G_2$ [Sh89].  Here is our main result:

\newtheorem*{Main Theorem}{Theorem}

\begin{Main Theorem} Let $\pi_1$ and $\pi_2$ be supercuspidal representations of $\GL_2(k)$ having the same central character.  Then for all sufficiently highly ramified characters $\omega$ of $\GL_2(k)$, we have
\[ \gamma(s,\pi_1 \otimes \omega, \operatorname{Sym}^3,\psi) = \gamma(s,\pi_2 \otimes \omega, \operatorname{Sym}^3,\psi). \]

\end{Main Theorem}

Our main result is equivalent to Theorem 2.2.1, which is the stability of local coefficients for $\operatorname{GL}_2$ inside $G_2$.  This theorem is not a new result: it follows from the strong transfer of cusp forms on $\GL_2 \times \GL_3$ to automorphic forms on $\GL_6$ via Rankin products done by Kim and Shahidi in [KiSh00].  But our method of proof gives a purely local proof of this result, and by the general argument in [Sh12], yields a new proof of the equality of $\gamma(s,\pi, \operatorname{Sym}^3,\psi)$ with the corresponding Artin factor under the local Langlands correspondence.  

Our method of proof follows that of [CoShTs17].  We apply Shahidi's local coefficient formula to write $\gamma(s,\pi,\operatorname{Sym}^3,\psi)$ as a Mellin transform of a partial Bessel function, and then prove an asymptotic expansion of partial Bessel functions in our case.  We use a concrete realization of $G_2$ as a group of 7 by 7 matrices.  Our matrix computations were done using SAGE Math.  

We hope that our work here may give some insight into the general problem of proving stability through asymptotic expansion, especially for exceptional Lie groups.  A particularly difficult case we hope to tackle in the future is stability for the exterior cube gamma factor for $\operatorname{GL}_6$, which occurs in the exceptional Lie group $E_6$.

\section{Structure of the group $G_2$}

Throughout the paper, $k$ denotes a $p$-adic field.  In this first section, we realize the split form $\mathbf G$ of $G_2$ as a group of $7$ by $7$ matrices.  We then realize $\mathbf M = \operatorname{GL}_2$ as a maximal Levi subgroup of $\mathbf G$ corresponding to the short root.  It is this inclusion of Levi subgroups which allows us to define the symmetric cube gamma factor via the Langlands-Shahidi method.

We will make several constructions in $\mathbf G$, including Weyl group representatives and orbit space measures, which will used to apply Shahidi's local coefficient formula (Theorem 6.2 of [Sh02]) to calculate the symmetric cube gamma factor in a way that will be amenable to a proof of stability.  

If $\mathbf H$ is an algebraic group defined over $k$, then $H = \mathbf H(k)$ will denote the group of its rational points.

\subsection{Definition of the Lie algebra $G_2$}

Let $\mathfrak g$ be the split semisimple Lie algebra over $k$ of type $\mathrm{G}_2$.  We may realize $\mathfrak g$ as the set of $7$ by $7$ matrices of the form
\[ \begin{pmatrix} a & x_{01} & x_{32} & 0 & -x_{21} & y_{10} & 2x_{11} \\ y_{01} & b & x_{31} & x_{21} & 0 &  -y_{11} & 2x_{10} \\  y_{32} & y_{31} & (-a - b) & -y_{10} & y_{11}& 0 & 2y_{21} \\ 0 & y_{21} & -x_{10} & -a & -y_{01} & -y_{32} & 2y_{11} \\ -y_{21} & 0  & x_{11} &  -x_{01} & -b  & -y_{31} & 2y_{10} \\ x_{10} &  -x_{11} & 0 & -x_{32} & -x_{31} & (a+b) & 2x_{21} \\ y_{11} & y_{10} & x_{21} & x_{11} & x_{10} & y_{21} & 0 \end{pmatrix}. \tag{1.1.1} \]
A Cartan subalgebra of $\mathfrak g$ is 
\[ \mathfrak t = \{ t = \operatorname{diag}(a,b,-a-b,-a,-b,a+b,0)\}.\]
If $\alpha, \beta \in \mathfrak t^{\ast}$ are the elements in the dual of $\mathfrak t$ sending the $t$ above to $b, a - b$ respectively, then we see that the set of roots of $\mathfrak t$ in $\mathfrak g$ is given by $\pm \Phi^+$, where
\[ \Phi^+ = \{ \alpha, \beta, \alpha + \beta, 2\alpha + \beta, 3\alpha + \beta, 3\alpha + 2\beta \} \]
is a set of positive roots, with simple roots $\Delta = \{\alpha, \beta \}$.

The positive eigenspaces are given by the variables $x_{ij}$, and the negative eigenspaces are given by the $y_{ij}$; for example, $x_{32}$ corresponds to the root $3\alpha+2\beta$, and $y_{11}$ corresponds to the root $-(\alpha+\beta)$.  

The variables $x_{ij}$ also fix positive root vectors $\mathbf x_{ij}: \mathbb G_a \rightarrow \mathbf U_{i \alpha + j \beta}$ in a natural way: $\mathbf x_{ij}(t)$ sends $t$ to the exponential of the matrix in (1.1) with $x_{ij} = t$ and all other entries zero.

\subsection{A matrix realization of $\mathbf G$}

Now that we have realized $\mathfrak g$ as a Lie algebra of matrices, we can define $\mathbf G$ to be the connected algebraic group in $\operatorname{GL}_{7,k}$ with Lie algebra $\mathfrak g$.  Let $\mathbf T$ be the maximal torus of $\mathbf G$ whose Lie algebra is $\mathfrak t$, and let $\mathbf B = \mathbf T \mathbf U$ be the Borel subgroup of $G$ corresponding to the system of positive roots in (1.1).  

Let $\mathbf P = \mathbf M \mathbf N$ be the parabolic subgroup of $G$ corresponding to the simple root $\beta$.  The Levi subgroup $\mathbf M$ is the centralizer of the kernel of $\beta$, regarded as a rational character of $ \mathbf T$.  We see that $\mathbf M$ is isomorphic to $ \GL_2$, and consists of matrices of the form
\[ \begin{pmatrix} A \\  & \Det A^{-1} \\ & & ^tA^{-1} \\ & & & \Det A \\ & & & &1 \end{pmatrix} \tag{1.2.1} \]
where $A \in \operatorname{GL}_2$.  The unipotent radical $\mathbf N$ of $\mathbf P$ consists of matrices of the form $\exp(X)$, for $X \in \mathfrak g$ satisfying $y_{ij} = a = b = x_{01} = 0$.  It follows that elements of $P$ looks like
\[ \begin{pmatrix} \ast & \ast & \ast & \ast & \ast & 0 & \ast \\ \ast & \ast & \ast & \ast & 0 & 0 & \ast \\ 0 & 0 & \ast & 0 & 0 & 0 & 0 \\ 0 & 0 & \ast & \ast & \ast & 0 & 0 \\ 0 & 0 & \ast & \ast & \ast & 0 &  0 \\ \ast & \ast & \ast & \ast & \ast & \ast & \ast \\ 0 & 0 & \ast & \ast & \ast & 0 & 1\end{pmatrix}. \tag{1.2.2} \]

\subsection{Splitting and Weyl group representatives}

Recall that $T = \mathbf T(k), \mathbf G = G(k)$ etc.  For each root $\gamma$ of $T$ in $G$, let $\mathbf U_{\gamma}$ be the corresponding root subgroup.  The variables $x_{ij}$ in (1.1) define root vectors $\mathbf x_{\gamma}: \mathbb G_a \rightarrow \mathbf U_{\gamma}$ for each positive root $\gamma$.  The choice of simple root vectors define canonical representatives $\dot w \in N_G(T)$ of each element $w$ in the Weyl group $W = N_G(T)/T$.  

This is done as follows:     For each simple root $\gamma \in \Delta$, there is a unique root vector $\mathbf x_{-\gamma}: \mathbb G_a \rightarrow U_{-\gamma}$ such that $\dot w_{\gamma} = \mathbf x_{\gamma}(1)\mathbf x_{-\gamma}(1)\mathbf x_{\gamma}(1)$ lies in the normalizer of $T$.  Then $\dot w_{\gamma}$ will be the canonical representative of $w_{\gamma}$.  We easily compute:
\[ \dot w_{\alpha} = \begin{pmatrix} & & & & & -1 \\ & & & & 1 \\ & & & 1 \\ & & -1 \\ & 1 \\ 1 \\ & & & & & & -1 \end{pmatrix} \]
\[ \dot w_{\beta} = \begin{pmatrix} 0 & 1 \\ -1 & 0 \\ & &  1 \\ &  & & 0& 1 \\ &  & & -1 & 0 \\ & & & & & 1 \\ & & & & & & 1\end{pmatrix}.\]
For a non-simple reflection $w$, the representative $\dot w$ is defined in terms of a reduced decomposition of $w$: if $(w_1, ... , w_r)$ is a reduced decomposition of $w$, where $w_i$ is equal to either $w_{\alpha}$ or $w_{\beta}$, then we set $\dot w = \dot w_1 \cdots \dot w_r$.  This will be independent of the choice of reduced decomposition. In particular, the long element $w_l$ has reduced decomposition $w_{\alpha}w_{\beta}w_{\alpha}w_{\beta}w_{\alpha}w_{\beta}$, so we have
\[ \dot w_l = \begin{pmatrix} & & & 1 \\ & & & & 1 \\ & & & & & 1\\ 1 \\ & 1 \\ & & 1 \\ & & & & & & -1 \end{pmatrix}\]
Finally, let $w_0$ be the unique element of $W$ which sends $\beta$ to a simple root and $\alpha$ to a negative root.  Explicitly, $w_0 = w_l w_{\beta}$, and in fact $w_0(\beta) = \beta$.  We have
\[ \dot w_0 = \dot w_l \dot w_{\beta}^{-1} = \begin{pmatrix} & & & & -1 \\ & & & 1 \\ & & & & & 1 \\ & -1 \\ 1 \\ & &  1 \\ & & & & & & -1 \end{pmatrix}.\]
The choice of splitting also defines a \emph{generic character} of $U$.  Let $\psi$ be a fixed nontrivial character of $k$.  If $u \in U$, then $u = \exp(\mathbf u)$ for a unique $\mathbf u \in \mathfrak g$.  In (1.1.1), the entries of $\mathbf u$ are zero except for possibly the $x_{ij}$.  We use $\psi$ to define a character of $U$ by the formula
\[ u \mapsto \psi(x_{01} + x_{10}). \]
We will also call this character $\psi$.  

\subsection{The normalized unramified character}

Let $\rho$ be half the sum of the roots of $\mathbf T$ in $\mathbf N$.  We have $2 \rho = 10\alpha + 5\beta$.  Let $X(\mathbf T)$ be the group of rational characters of $\mathbf T$.  It has $\alpha, \beta$ as a basis.  The Weyl group $W = N_G(T)/T$ acts on $X(\mathbf T)$ by 
\[ w.\gamma(t) = \gamma(n^{-1}tn)\]
where $n \in N_G(T)$ represents $w$.  We have that 
\[ w_{\alpha}(\beta) = 3\alpha+\beta \tag{1.4.1} \]
\[ w_{\beta}(\alpha) = \alpha+\beta \tag{1.4.2}\]
We define a symmetric, positive definite bilinear form on $X(\mathbf T) \otimes_{\mathbb Z} \mathbb R$ by
\[ (a \alpha + b \beta, a'\alpha + b'\beta) = aa' + 3bb' - \frac{3}{2}ab' - \frac{3}{2}a'b. \]
On account of (1.4.1) and (1.4.2), this form is easily seen to be invariant under the action of $W$.  

As in [Sh02], we set
\[ \tilde{\alpha} = \langle \rho, \alpha \rangle^{-1}\rho = 2 \frac{(\rho,\alpha)}{(\alpha,\alpha)}\rho = 4\rho. \]
Then for $s \in \mathbb C$, we get an unramified character of $M(k) = \GL_2(k)$, defined by
\[ m \mapsto q^{\langle s \tilde{\alpha}, H_M(m) \rangle} = |\det(m)|^{10s}.\]

\subsection{Orbit space representative}

Let $\mathbf U_{\mathbf M} = \mathbf U \cap \mathbf M$.  The Lie algebra of $\mathbf U_{\mathbf M}$ is one-dimensional and consists of all matrices of the form
\[ \mathbf u =  \begin{pmatrix} 0 & x \\ 0 & 0 \\ & & 0 \\ & & & 0 & 0 \\ & & & -x & 0 \\ & & & & & 0 \\ & & & & & & 0 \end{pmatrix}. \]
We consider the action of $U_M$ on $N$ by conjugation.  We are interested in finding a suitable open dense subset $N'$ of $N$, stable under the action of $U_M$, such that the quotient space $U_M \backslash N'$ has a nice $p$-adic manifold structure.  We will also be interested in finding a measure $d \dot n$ on $U_M \backslash N'$, such that integration over $N$ can be recovered by double integration over $U_M$ and $U_M \backslash N'$.  

Let us consider a typical element $n = \exp(\mathbf n)$ of $N$, for
\[ \mathbf n = \begin{pmatrix} 0 & x_{01} & x_{32} & 0 & -x_{21} & 0 & 2x_{11} \\ 0 & 0 & x_{31} & x_{21} & 0 & 0 & 2x_{10} \\  0 & 0 &0 & 0 & 0& 0 & 0 \\ 0 & 0 & -x_{10} &0 &0 & 0 & 0 \\ 0 & 0  & x_{11} &  -x_{01} & 0  & 0 & 0 \\ x_{10} &  -x_{11} & 0 & -x_{32} & -x_{31} & 0 & 2x_{21} \\ 0 & 0 & x_{21} & x_{11} & x_{10} & 0 & 0 \end{pmatrix} \]
To simplify the notation, we can identify $\mathbf n$ with a tuple
\[ \mathbf n = (x_{01}, x_{11}, x_{21}, x_{31}, x_{32}) \tag{1.5.1}\]
If $u = \exp(\mathbf u)$, then we have
\[ unu^{-1} = \exp(\mathbf n_1)\]
where
\[ \mathbf n_1 = (x_{10}, xx_{10}+x_{11}, x_{21}, x_{31}, xx_{31} + x_{32}) \]
We see that on the open dense subset $N'$ of $N$ consisting of elements $\exp(\mathbf n)$ such that $x_{10} \neq 0$, the action of $U_M$ is simple.  Taking $x = -x_{11}x_{10}^{-1}$, we get that every element of $N'$ is conjugate by a unique element of $U_M$ to a unique element of the form $\exp(\mathbf n_0)$, for
\[ \mathbf n_0 = (x_{10}, 0, x_{21}, x_{31}, x_{32}) \]

Putting this together, we have:

\newtheorem{1_5_1}{Lemma}[subsection]

\begin{1_5_1} Let $D$ be the set of $\exp(\mathbf n_0) \in N$, for
\[ \mathbf n_0 = (x_{10},0,x_{21},x_{31},x_{32})\]
Then $U_M \times D \rightarrow N', (u,n) \mapsto unu^{-1}$ is an isomorphism of analytic manifolds.  In particular, the map sending $n \in N'$ to its unique conjugate in $D$ is a submersion of manifolds, so $D$ is the quotient of $N'$ under the action of $U_M$ in the category of analytic manifolds.

\end{1_5_1}

Next, let $f$ be a smooth, compactly supported function on $N$.  We claim there is a measure $d \dot n$ on $D$ such that integration over $N$ can be recovered by double integration over $U_M$ and $D$.  

We take $d \dot n$ to be the measure induced from the exponential map and the measure $|x_{10}|dx_1dx_3dx_4dx_5$ on $\log D$.  It suffices to compute everything at the level of the Lie algebra: a typical element of $\log D$ can be written as $\mathbf n_0 = (x_{10},0,x_{21},x_{31},x_{32})$.  If $u = \exp(\mathbf n) \in U_M$, then 

\[ u \mathbf n_0 u^{-1} = (x_{10},xx_{10},x_{21},x_{31},xx_{31}+x_{32})\]
so
\begin{equation*}
    \begin{split}
        \int\limits_{U_M} \int\limits_{\log D} f(u \mathbf n_0 u^{-1})d(\log \dot n) du  = \int\limits_k \int\limits_{k^4} f((x_{10},xx_{10},x_{21},x_{31},xx_{31}+x_{32}))\\|x_{10}|dx_1dx_3dx_4dx_5 dx.
    \end{split}
\end{equation*} 
We change $x_{32}$ to $x_{32} - xx_{31}$, and then change $x$ to $xx_{10}^{-1}$.  This becomes
\[ \int\limits_{k^5} f(x_{10},x,x_{21},x_{31},x_{32})dx_1dx_3dx_4dx_5 dx = \int\limits_{\mathfrak n} f(n)dn.\]

We have shown:

\newtheorem{1_5_2}[1_5_1]{Lemma}

\begin{1_5_2} Let $D$ be as in Lemma 1.5.1, and let $d \dot n$ be the measure on $D$ induced from the exponential map and the measure $|x_{10}| dx_{10} dx_{21} dx_{31} dx_{32}$ on $\log(D)$.  Then
\[ \int\limits_N f(n) dn = \int\limits_{U_M} \int\limits_D f(un_0u^{-1}) d \dot n_0 du\]
for all $f \in \mathscr C_c^{\infty}(N)$.

\end{1_5_2}

We next introduce the conjugation action of the center $Z_M$ of $M$.  By (1.2.1), $Z_M$ consists of all matrices of the form
\[ z = \operatorname{diag}(t,t,t^{-2},t^{-1},t^{-1},t^2,1)\]
for $t \in k^{\ast}$, and we can therefore identify $Z_M$ with $k^{\ast}$.  If $\mathbf n = (x_{10},x_{11},x_{21},x_{31},x_{32})$ is an element of the Lie algebra of $N$ (as in equation (1.2.1)), we have
\[ z \mathbf n z^{-1} = (tx_{10},tx_{11},t^2x_{21},t^3x_{31},t^3x_{32}). \]
The action of $Z_M$ on $N$ commutes with that of $U_M$, and $U_M \times Z_M$ acts simply on $N'$.  Similar to Lemmas 1.5.1 and 1.5.2, we can further quotient out $D$ by $Z_M$ to obtain a fundamental domain $D_0$ for $U_M \times Z_M$ and, with an appropriate measure on $D_0$, we may recover integration over $N$ by integration of $U_M \times Z_M \times D_0$.  The proof is similar to that of the previous two lemmas, so we omit the details.

\newtheorem{1_5_3}[1_5_1]{Lemma}

\begin{1_5_3} Let $D_0$ be the set of $\exp(\mathbf n_0) \in N$, for
\[ \mathbf n_0 = (1,0,x_{21},x_{31},x_{32})\]
(i): The map $U_M \times Z_M \times D_0 \rightarrow N', (u,z,n) \mapsto uznz^{-1}u^{-1}$ is an isomorphism of analytic manifolds.  In particular, the map sending $n \in N'$ to its unique representative in $D_0$ is a submersion of manifolds, so $D_0$ is the quotient of $N'$ by $Z_M \times U_M$ in the category of analytic manifolds.  

(ii): With the measure $d \dot n = dx_{21} dx_{31} dx_{32}$ on $D_0$, integration over $N$ can be recovered by integration over $U_M \times Z_M \times D_0$.  That is, 
\[ \int\limits_N f(n) dn = \int\limits_{Z_M} \int\limits_{U_M} \int\limits_{D_0} f(zunu^{-1}z^{-1})q^{\langle 2\rho,H_M(z)\rangle} d \dot n du \]
for all $f \in \mathscr C_c^{\infty}(N)$.  
\end{1_5_3}

Note that if $z$ identifies with $t \in k^{\ast}$, then $q^{\langle 2\rho, H_M(z) \rangle} = |t|^{10}$ by (1.4).  

Since the preimage of dense open sets under a submersion of manifolds remain dense open, it is clear that Lemma 1.5.3 holds if $D_0$ is replaced by any open dense subset $W$ of $D_0$, and $N'$ is replaced by the orbit of $W$ under $U_M \times Z_M$. In the notation of [Sh02], $D_0$ would be written as $Z_M U_M \backslash N$.  

\subsection{The decomposition $\dot w_0^{-1}n = mn'\bar{n}$}

Let $\overline{\mathbf N}$ be the unipotent radical of the parabolic subgroup opposite to $\mathbf P$.  The set of $n \in N$ such that $\dot w_0^{-1}n \in P \overline{N}$ is open and dense in $N$.  We may then write $\dot w_0^{-1}n$ uniquely as $mn'\bar{n}$ for $m \in M, n \in N, \bar{n} \in \overline{N}$.  

Let us consider $n$ in the fundamental domain $D_0$ of (1.5) under the action of $Z_M U_M$, so $n = \exp(\mathbf n)$, for $\mathbf n = (1,0,x_{21},x_{31},x_{32})$.

Let $\bar{n} \in \overline{N}$, for $\bar{n} = \exp(\bar{\mathbf n})$, and
\[ \bar{\mathbf n} =  \begin{pmatrix} 0 & 0 & 0 & 0 & 0 & y_{10} & 0 \\ 0 & 0 & 0 & 0 & 0 &  -y_{11} & 0 \\  y_{32} & y_{31} & 0 & -y_{10} & y_{11}& 0 & 2y_{21} \\ 0 & y_{21} & 0 & 0 & 0 & -y_{32} & 2y_{11} \\ -y_{21} & 0  & 0 &  0 & 0  & -y_{31} & 2y_{10} \\ 0 &  0 & 0 & 0 & 0 & 0 & 0 \\ y_{11} & y_{10} & 0 & 0 & 0 & y_{21} & 0 \end{pmatrix}.  \]
Under the dense open condition $x_{21} \neq 0, x_{32} \neq 0$ and $D := x_{21}^2 + x_{32} \neq 0$, we see that $\dot w_0^{-1}n$ does lie in $P \overline{N}$.

Indeed, if we set
\[ y_{10} = \frac{-x_{32}}{D} \]
\[ y_{11} = \frac{-x_{31} + \frac{1}{2}x_{21}}{D} \]
\[ y_{21} = \frac{-x_{21}}{D}\]
\[ y_{31} = \frac{\frac{1}{2}x_{32}x_{21}}{D^2}\]
\[ y_{32} = \frac{\frac{3}{4}x_{21}^2+x_{10}^2x_{32}+\frac{1}{2}x_{21}x_{31}}{D^2}\]
then we have $p = \dot w_0^{-1}n \bar{n}^{-1} \in P$.  This can be seen by looking at which entries a typical element of $P$ must be zero.  Writing $p = mn$ with $m \in \GL_2(k)$, the entries of $m$ can be read off directly from the four upper left entries of $p$:
\[ m = \begin{pmatrix} a & b \\ c & d \end{pmatrix} \]
\[ a = y_{10}y_{11} + y_{21}\]
\[ b = y_{10}^2\]
\[ c = -y_{11}^2 + \frac{1}{2}y_{11}y_{21} - y_{32}\]
\[ d = -y_{10}y_{11} + \frac{1}{2}y_{10}y_{11} + y_{21} - y_{31}.\]
If we write $m$ in the Bruhat decomposition $m = u_1 \begin{pmatrix} 0 & 1 \\ -1 & 0 \end{pmatrix} t u_2$, then
\[ u_1 = \begin{pmatrix} 1 & \frac{a}{c} \\& 1 \end{pmatrix}\]
\[ u_2 = \begin{pmatrix} 1 & \frac{d}{c}\\ & 1 \end{pmatrix} \]
\[ t = \begin{pmatrix} - \frac{1}{c} \det(m) \\ & -c \end{pmatrix}. \]
The following result comes from the unexpected observation that the entries of $m, u_1, u_2,$ and $t$ are quotients of homogeneous polynomials in the variables $x_{21}, x_{31},$ and $x_{32}$.  We do not currently have a good explanation for why this homogeneity occurs in $G_2$, but we expect it to occur in other Lie groups.

\newtheorem{1_5_4}[1_5_1]{Lemma}

\begin{1_5_4} If $t \in k^{\ast}$, and $(x_{21}, x_{31}, x_{32})$ is changed to $(t x_{21},tx_{31},t^2x_{32})$, then $y_{10}, y_{11}, y_{21},y_{31},y_{32}$ will change respectively to $y_{10}, t^{-1}y_{11}, t^{-1}y_{21}, t^{-1}y_{31}, t^{-2}y_{32}$.  The matrix $t = \operatorname{diag}(t_1,t_2)$ will change to $\operatorname{diag}(t_1, t^{-2}t_2)$.  If we write \[u_1 = \begin{pmatrix} 1 & x_1 \\ & 1 \end{pmatrix} \, \, \, u_2 = \begin{pmatrix} 1 & x_2 \\ & 1 \end{pmatrix}\] then these change respectively to 
\[ \begin{pmatrix} 1 & tx_1 \\ & 1 \end{pmatrix} \, \, \, \,  \, \begin{pmatrix} 1 & tx_2 \\ & \end{pmatrix} \]

\end{1_5_4}

\begin{proof} Direct computation.

\end{proof}
\subsection{Open compact subgroups}

We will need a nice collection of open compact subgroups $\overline{N}_{\kappa} : \kappa \geq 1$ of $\overline{N}$, such that $z \overline{N}_{\kappa} z^{-1}$ only depends on $|t|$ for $z = \begin{pmatrix} t \\& t \end{pmatrix} \in Z_M$.

Let us consider two arbitrary elements $\bar{n} = \exp(\mathbf n)$ and $\bar{n}' = \exp(\mathbf n')$ of $\overline{N}$ for $\mathbf n, \mathbf n'$ in the Lie algebra of $\overline{N}$.  If we write
\[ \mathbf n = (y_{10}, y_{11}, y_{21}, y_{31}, y_{32}) \]
\[ \mathbf n' = (y_{10}', y_{11}', y_{21}', y_{31}', y_{32}') \]
then we see that 
\[ \bar{n} \bar{n}' = \exp(z_{10}, z_{11}, z_{21}, z_{31}, z_{32}) \]
where \[z_{10} = y_{10} + y_{10}'\] \[z_{11} = y_{11} + y_{11}'\] \[z_{21} = y_{11}y_{10}' - y_{10} y_{11}' + y_{21} + y_{21}'\]
\[ z_{31} = \frac{1}{2}(-y_{10} y_{11}z_{10}' +  y_{11} y_{10}'^2 + y_{10}^2y_{11}' - y_{10}y_{10}'y_{11}') + \frac{3}{2}(y_{21} y_{10}' - y_{10}y_{21}') + y_{31} + y_{31}' \]
\[ z_{32} = \frac{1}{2}(-y_{11}^2y_{10}' + y_{10} y_{11}y_{11}' + y_{11} y_{10}'y_{11}' - y_{10}y_{11}'^2 - y_{11}y_{21}) + \frac{3}{2} y_{21}y_{11}' - y_{11}y_{21}' + y_{32} + y_{32}'.\]

For $\kappa \geq 1$, we will define $\overline{N}_{\kappa}$ to be the set of $\exp(\mathbf n)$ for
\[ \mathbf n \in (\mathfrak p_k^{-\kappa^2}, \mathfrak p_k^{-\kappa}, \mathfrak p_k^{-\kappa^3}, \mathfrak p_k^{-\kappa^5}, \mathfrak p_k^{-\kappa^4}). \]
On account of the potentially large discrepancies between the absolute values of the coordinates of entires of $\overline{N}_{\kappa}$, and the multiplication formula in $\overline{N}$ given above, we can deduce the following lemma.

\newtheorem{1_7_1}{Lemma}[subsection]

\begin{1_7_1} (i): For large enough $\kappa$, the sets $\overline{N}_{\kappa}$ are open compact subgroups of $\overline{N}$, such that for \[z = \begin{pmatrix} t \\ & t \end{pmatrix} \in Z_M,\] the subgroup $z \overline{N}_{\kappa}z^{-1}$ only depends on $|t|$.

(ii): Let $U_1$ be an open compact subgroup of $U_M$, and let $\varphi_{\kappa}$ be the characteristic function of $\overline{N}_{\kappa}$.  There exists a $\kappa_0$ such that for all $\kappa \geq \kappa_0$, we have
\[ \varphi_{\kappa}(u\bar{n}u^{-1}) = \varphi_{\kappa}(\bar{n})\]
for all $\bar{n} \in \overline{N}$ and all $u \in U_1$.

\end{1_7_1}
\section{Partial Bessel functions}

In this section, we apply Shahidi's local coefficient formula to calculate the symmetric cube gamma factor as a Mellin transform of a partial Bessel function, up to some Tate gamma factors.  We then develop the asymptotics of this partial Bessel function to prove the stability result (Theorem 2.2.1).

\subsection{Definition of the partial Bessel function}

Let $\omega: Z_M \rightarrow \mathbb C^{\ast}$ be any character, and let $f: M \rightarrow \mathbb C$ be any locally constant function which is compactly supported modulo $Z_M$ and which transforms according to $\omega$, that is $f(zg) = \omega(z) f(g)$ for all $z \in Z_M$ and $g \in G$.  We denote the space of such functions by $\mathscr C_c^{\infty}(M ; \omega)$.  For example, $f$ could be a matrix coefficient of a supercuspidal representation of $M = \GL_2(k)$.  We define $W^f: M \rightarrow \mathbb C$ by
\[ W^f(m) = \int\limits_{U_M} f(xm) \overline{\psi(x)} dx.\]
Since $f$ is compactly supported modulo $Z$, this integral converges absolutely.  Let $n \in N$ be an element for which $\dot w_0^{-1}n \in P\overline{N}$.  Writing $\dot w_0^{-1}n = mn'\bar{n}$ as in (1.6), and letting $\varphi_X$ be the characteristic function of an open compact set $X$ in $\overline{N}$, we define the \emph{partial Bessel function} 
\[ J_{\varphi_X}(n,f) = \int\limits_{U_M} W^f(mu) \varphi_X(u \bar{n} u^{-1}) \overline{\psi(u)} du. \]
This integral converges absolutely, because $u \mapsto u \bar{n} u^{-1}$ is a homeomorphism of $U_M$ onto the orbit of $\bar{n}$ under the conjugation action of $U_M$ on $\overline{N}$, and orbits of unipotent algebraic groups on affine varieties are closed.

When $f$ is supported inside the big cell, a smoothness property holds for the arguments of partial Bessel function.

\newtheorem{2_1_1}{Proposition}[subsection]

\begin{2_1_1} Assume that $f$ vanishes on $B_M = U_MT$.  There exists an open compact subgroup $H$ of $k^{\ast}$ depending on $f$ and $\omega$, such that the following holds: if $t \in H$, $n = \exp (1,0, x_{11}, x_{31}, x_{32})$, and $n_1 = \exp(1,0, tx_{11}, tx_{31}, t^2x_{32})$, then 
\[ J_{\varphi_{\kappa}}(n,f) = J_{\varphi_{\kappa}}(n',f)\]
for all $\kappa$.  

\end{2_1_1}

\begin{proof} Our hypothesis is that $f$ is supported inside the \emph{big cell} of $M$, that is $f(b) = 0$ for all $b \in B_M$.  Then $f$ is compactly supported modulo $Z_M$ as a function on the big cell $U_M\dot w_M TU_M$.  Let $T' = \{ \operatorname{diag}(1,t_2) : t_2 \in k^{\ast} \} \subset T$, so that $T$ is the direct product of $Z_M$ and $T'$.  

There exist open compact subgroups $U_1$ and $U_2$ of $U_M$, and a compact set $\Omega$ of $T'$, such that if $f(u_1 \dot w_M zt'u_2) \neq 0$ for $u_i \in U_M, z \in Z_M,$ and $t' \in T'$, then $u_i \in U_i$ and $t' \in \Omega$.  It follows from this uniformity that there must be a small neighborhood $H$ of $1$ in $k^{\ast}$, which we may take to be a compact open subgroup, such that if $n$ is changed to $n_1$ by $t \in H$, then for the corresponding decomposition $\dot w_0^{-1}n_1 = m_1 n_1' \bar{n}_1'$, we have $f(xmu) = f(xm_1u)$ for all $x, u \in U_M$.  This follows from the fact that $f$ is locally constant and $m,n', \bar{n}$ are continuous functions of $n$, or more explicitly from Lemma 1.6.1.

We also see that since $t$ is necessarily in $\mathcal O_k^{\ast}$, this will not affect the calculation of the characteristic function $\varphi_{\kappa}$, that is we will have $\varphi_{\kappa}(u\bar{n}u^{-1}) = \varphi_{\kappa}(u \bar{n}'u^{-1})$ for all $u \in U_M$.  \end{proof}

\subsection{Local coefficient formula}

Let $\pi$ be a generic, irreducible representation of $M = \GL_2(k)$.  Let $C_{\psi}(s,\pi)$ be the Shahidi local coefficient attached to $\pi$, relative to $M$ inside $G$ ([Sh90], [Sh02]).  The symmetric cube gamma factor $\gamma(s,\pi, \operatorname{Sym}^3,\psi)$, as defined by the Langlands-Shahidi method, is related to the local coefficient by the formula
\[ C_{\psi}(s,\pi \otimes \omega_{\pi}) = \gamma(2s, \omega_{\pi}^3,\psi) \gamma(s,\pi, \operatorname{Sym}^3,\psi). \tag{2.2.1}\](Corollary to Proposition 2.2 of [Sh89]).  Here $\omega_{\pi}$ is the central character of $\pi$, and $\gamma(s, \omega_{\pi}^3,\psi)$ is the Tate gamma factor attached to the character $\omega_{\pi}^3$ of $k^{\ast}$.  Our main theorem is therefore equivalent to the stability of local coefficients:

\newtheorem{2_2_1}{Theorem}[subsection]

\begin{2_2_1} Let $\pi_1$ and $\pi_2$ be supercuspidal representations of $\GL_2(k)$ with the same central character.  Then for all sufficiently highly ramified characters $\omega$ of $\GL_2(k)$, we have
\[ C_{\psi}(s, \pi_1 \otimes \omega) = C_{\psi}(s, \pi_2 \otimes \omega). \]

\end{2_2_1}

Let $\pi$ be a supercuspidal representation of $\GL_2(k)$ with ramified central character $\omega_{\pi}$.  Theorem 6.2 of [Sh02] gives us the following formula for $C_{\psi}(s,\pi)^{-1}$ as a Mellin transform of a partial Bessel function:
\begin{equation*}
    \begin{split}
        C_{\psi}(s,\pi)^{-1} = \gamma(2 \langle \tilde{\alpha},\alpha^{\vee} \rangle s, \omega_{\pi}(\dot w_0 \omega_{\pi}^{-1}),\psi)^{-1} \\ \int\limits_{Z_M U_M \backslash N} J_{\varphi_{z_0\overline{N}_0z_0^{-1}}}(n,f) \omega_{\pi}^{-1}( \dot w_0 \omega_{\pi})(x_{\alpha}) q^{\langle s \tilde{\alpha} + \rho, H_M(m) \rangle} d \dot n
    \end{split}
\end{equation*}
We will explain the notation in this formula before simplifying it.  Here $\alpha^{\vee}$ is a coroot of $Z_M$ with the property that $\langle \alpha, \alpha^{\vee} \rangle = 1$.  By (1.4), the pairing $\langle \tilde{\alpha},\alpha \rangle$ equals $20$.  The expression $\omega_{\pi}(\dot w_0 \omega_{\pi}^{-1})$ denotes the character $z \mapsto \omega_{\pi}(z) \omega_{\pi}^{-1}(\dot w_0^{-1}z\dot w_0)$ of $Z_M$, which is this case is equal to just $\omega_{\pi}^2$.  

The integral is over the quotient of an open dense subset of $N$ under the conjugation action of $Z_MU_M$, and the measure $d\dot n$ is the ``orbit space measure'' which allows us to recover integration over $N$ by integration over $Z_M, U_M$, and $Z_M U_M \backslash N$.  We have already identified the space $Z_M U_M \backslash N$ and the measure $d \dot n$ in Lemma 1.5.3: we can take $Z_M U_M \backslash N$ to be the torus $R = \{ \exp(1, 0,x_{21}, x_{31}, x_{32}) : x_{ij} \neq 0\}$, which is actually a subset of $N$.  The measure $d \dot n$ is given by \[ d \dot n = dx_{21} d x_{31} d x_{32}. \]
For $n \in R$, we can write $\dot w_0^{-1}n = mn' \bar{n}$ whenever the ``discriminant'' $D = x_{21}^2 + x_{32}$ is not zero.  The set of $n \in R$ for which $D  = 0$ is of measure zero, and we can declare the integrand to just be $0$  for such $n$.  

In the partial Bessel function $J_{\varphi_{z\overline{N}_0z^{-1}}}(n,f)$, $f$ is a matrix coefficient of $\pi$ with $f(e) = 1$, $z_0 = t_0 I_2$ is a diagonal matrix in $Z_M$, where $t_0$ is an element of $k^{\ast}$ whose absolute value depends on the conductors of  $\omega_{\pi}$ and $\psi$, and $\overline{N}_0$ is an open compact subgroup of $\overline{N}$ with the property that $z \overline{N}z^{-1}$ depends only on $|t|$ for all $z = t I_2 \in Z_M$.  The open compact subgroup $\overline{N}_0$ depends on $\pi$, but once found, may be replaced by any larger open compact subgroup of $\overline{N}$ with the same property.  

Finally, if $\dot w_0^{-1}n  = mn' \bar{n}$, then $\dot w_0^{-1} \bar{n} \dot w_0 \in N$, which we may write as \[\dot w_0^{-1} \bar{n} \dot w_0 = \exp(c_{10},c_{11},c_{21},c_{31},c_{32})\] for $c_{ij} \in k$.  Then $x_{\alpha}$ designates $c_{10}$, which we calculate as
\[ x_{\alpha} = \frac{\frac{1}{2}x_{21}-x_{31}}{x_{21}^2+x_{32}}\]
Let $\{\overline{N}_{\kappa} \}$ be the open compact subgroups of (1.7).  Let $\kappa_0$ be a sufficiently large integer so that $z_0^{-1} \overline{N}_{\kappa_0} z_0$ contains $\overline{N}_0$.  Then for all $\kappa \geq \kappa_0$, we may use $\overline{N}_{\kappa}$ in place of $z_0 \overline{N}_0 z_0^{-1}$.  We arrive at the following reformulation of Theorem 6.2 of [Sh02] in our case:

\newtheorem{2_2_2}[2_2_1]{Proposition}

\begin{2_2_2} Let $\pi$ be an irreducible, supercuspidal representation of $\GL_2(k)$ with ramified central character.  Let $f$ be a matrix coefficient of $\pi$ with $W^f(e) = 1$.  Then there exists an integer $\kappa_0$ depending on $\pi$ such that for all $\kappa \geq \kappa_0$,
\begin{equation*}
    \begin{split}
         C_{\psi}(s,\pi)^{-1} &  = \gamma(40s, \omega_{\pi}^2,\psi)^{-1} \int\limits_R J_{\varphi_{\kappa}}(n,f) \omega_{\pi}^{-2}(\frac{\frac{1}{2}x_{21}-x_{31}}{x_{21}^2+x_{32}})\\ & |\det(m)|^{10s + \frac{5}{2}} dx_{21} dx_{31} dx_{32}.
    \end{split}
\end{equation*}

\end{2_2_2}

\subsection{Moving up to the big cell}

The smoothness property of Proposition 2.1.1 is crucial for our stability result.  It only holds for functions supported inside the big cell.  Our matrix coefficient $f$ occurring in the local coefficient formula of Proposition 2.2.2 are not supported in the big cell, since they are assumed to satisfy $W^f(e) = 1$.  

In order to access the smoothness result of Proposition 2.1.1., we will need to prove an asymptotic expansion formula of $J_{\varphi_{\kappa}}(n,f)$.  We are looking for two functions $f_1$ and $f_1^0$ in $\mathscr C_c^{\infty}(M ; \omega_{\pi})$ such that
\[ J_{\varphi_{\kappa}}(n,f) = J_{\varphi_{\kappa}}(n,f_1) + J_{\varphi_{\kappa}}(n,f_1^0).\]The function $f_1$ will only depend on $\omega_{\pi}$, and the second $f_1^0$ will be supported inside the big cell.  

The process of obtaining $f_1^0$ from $f$ and $f_1$ is done in this section.  It is very similar to $\S$ 5.4 of [CoShTs17].  

\newtheorem{3_1_1}{Lemma}[subsection]

\begin{3_1_1} Let $f \in \mathscr C_c^{\infty}(M ; \omega)$.  If $U_1$ and $U_2$ are open compact subgroups of $U_M = k$, define $f' \in \mathscr C_c^{\infty}(M ;\omega)$ by
\[ f'(m) = \frac{1}{\operatorname{meas}(U_1) \operatorname{meas}(U_2)} \int\limits_{U_1 \times U_2} f(u_1mu_2) \overline{\psi(u_1u_2)} du_1 du_2.\]
Then there exists a $\kappa_0$ depending on $U_2$ such that for all $\kappa \geq \kappa_0$, we have
\[ J_{\varphi_{\kappa}}(n,f) = J_{\varphi_{\kappa}}(n, f')\]
for all $n \in R$.  
\end{3_1_1}

\begin{proof} We take $\kappa$ sufficiently large so that $\varphi_{\kappa}(u_2 \bar{n} u_2^{-1}) = \varphi_{\kappa}(u)$ for all $u_2 \in U_2$ and $n \in R$ (Lemma 1.7.1).  Let $c = \frac{1}{\operatorname{meas}(U_1) \operatorname{meas}(U_2)}$, so that
\[ J_{\varphi}(n, f') = c  \int\limits_{U_M} \int\limits_{U_M} \int\limits_{U_1 \times U_2} f(xu_1mu_2 u) \varphi_{\kappa}(u \bar{n}u^{-1}) \overline{\psi(xu_1 u_2 u)} du_1 du_2 dx du.\]
We get the result by making the change of variables $x \mapsto xu_1^{-1}$ and $u \mapsto u_2^{-1}u$.
\end{proof}

We will use Lemma 2.3.1 to show that if the ``pure Bessel integral'' $W^f(-)$ vanishes on the center of $M$, then in calculating the partial Bessel integral $J_{\varphi}(-,f)$ we may replace $f$ by a function which is supported inside the big cell.  We do this in two steps (Lemma 2.3.2 and Proposition 2.3.3).

\newtheorem{3_1_2}[3_1_1]{Lemma}

\begin{3_1_2} Let $f \in \mathscr C_c^{\infty}(M ; \omega)$, and suppose that $W^f(e) = 0$.  Then there exists an $f_0 \in \mathscr C_c^{\infty}(M ;\omega)$ which vanishes on $U_M Z_M$, and an integer $\kappa_0$, such that 
\[ J_{\varphi_{\kappa}}(n, f) = J_{\varphi_{\kappa}}(n, f_0) \]
for all $n \in R$ and for all $\kappa \geq \kappa_0$.   \end{3_1_2}

\begin{proof} Consider the restriction of $f$ to $U_M T$.  Since $f$ is compactly supported modulo $Z_M$, there is an open compact subgroup $U_1$ of $U_M$ such that $f(ut) = 0$ implies $u \in U_1$.  If we set
\[ f_0(m) = \frac{1}{\operatorname{meas}(U_1)^2} \int\limits_{U_1 \times U_1} f(u_1m u_2) \overline{\psi(u_1u_2)} du_1 du_2 \]
then Lemma 2.3.1 tells us that for sufficiently large $\kappa$, we have $J_{\varphi_{\kappa}}(n,f) = J_{\varphi_{\kappa}}(n,f_0)$.  We need only show that $f_0(uz) = 0$ for all $z \in Z_M$ and $u \in U_M$.  If this is not the case, then there is a $z \in Z_M, u \in U_M$, and $u_1, u_2\in U_1$ such that $f(u_1uzu_2) \neq 0$.  This implies $u_1uu_2 \in U_1$, hence $u \in U_1$.  Now let $c = \operatorname{meas}(U_1)$, so that
\[ f_0(uz) = \frac{\omega(z)}{c^2} \int\limits_{U_1 \times U_1} f(u_1uu_2) \overline{\psi(u_1u_2)} du_1du_2.\]
We change variables $u_1 \mapsto u_1u^{-1}u_2^{-1}$, so that
\begin{equation*}
    \begin{split}
         f_0(uz) & = \frac{\omega(z)}{c^2} \psi(u) \int\limits_{U_1 \times U_1} f(u_1) \overline{\psi(u_1)} du_1du_2 \\
         & = \frac{\omega(z)}{c} \psi(u) \int\limits_{U_1} f(u_1) \overline{\psi(u_1)} du_1 \\
         & = \frac{\omega(z)}{c} \psi(u) W^f(e) \\
         & = 0.
    \end{split}
\end{equation*}
\end{proof}

\newtheorem{3_1_3}[3_1_1]{Proposition}

\begin{3_1_3} Let $f \in \mathscr C_c^{\infty}(M ; \omega)$, and suppose that $W^f(e) = 0$.  Then there exists an $f_0 \in \mathscr C_c^{\infty}(M ;\omega)$ which is supported inside the big cell, and an integer $\kappa_0$, such that 
\[ J_{\varphi_{\kappa}}(n, f) = J_{\varphi_{\kappa}}(n, f_0) \]
for all $n \in R$ and all $\kappa \geq \kappa_0$. 

\end{3_1_3}

\begin{proof} By Lemma 2.3.3, there is an $f_1 \in \mathscr C_c^{\infty}(M ; \omega)$ which vanishes on $U_MZ_M$ and satisfies $J_{\varphi}(n,f) = J_{\varphi}(n,f_1)$ for all $n \in R$.  We may therefore replace $f$ by $f_1$ and assume from the beginning that $f$ vanishes on $Z_M U_M$.

Identifying $M = \GL_2(k)$, we can write $T$ as the direct product of $Z_M$ and $T' = \{ \operatorname{diag}(1,y) : y \in k^{\ast} \}$.  Since $f$ is compactly supported modulo $Z_M$, there is an open subgroup $U_1$ of $U_M$, and a compact set $\Omega \subset T'$ such that if $f(ut') \neq 0$ for $u \in U_M$ and $t' \in T'$, then $u \in U_1$ and $t' \in \Omega$.  

Since we are assuming that $f$ vanishes on $U_MZ_M$, we can furthermore choose $\Omega$ to be disjoint from $Z_M$.  Identifying $T'$ with $k^{\ast}$ in the obvious way, we see that $\Omega$ is a compact set in $k^{\ast}$ which is bounded away from $1$.  Therefore each additive character $x \mapsto \psi(x(y^{-1}-1))$ for $y \in k^{\ast}$ is nontrivial, and we can find an open compact subgroup $U_2$ of $U_M$ such that
\[ \int\limits_{U_2} \overline{\psi(x(y-1))} dx = 0 \]
for all $y \in \Omega$.  We can also enlarge our subgroup $U_1$ so that $t' U_2 t'^{-1} \subset U_1$ for all $t' \in \Omega$.  Now set
\[ f_0(m) = \frac{1}{\operatorname{meas}(U_1) \operatorname{meas}(U_2)} \int\limits_{U_1 \times U_1} f(u_1mu_2) \overline{\psi(u_1 u_2)} du_1 du_2\]
so that for sufficiently large $\varphi$, we have $J_{\varphi}(n,f) = J_{\varphi}(n,f_0)$.  We are done if we can show that $f_0(ut) = 0$ for all $u \in U_M$ and $t \in T$.  Writing $t = zt'$ for $z \in Z_M$ and $t' \in T'$, we have $f_0(ut) = \omega(z) f_0(ut')$, so it suffices to show that $f_0(ut') = 0$.  Now $f_0(ut')$ is a scalar multiple of
\[ \int\limits_{U_1 \times U_2} f(u_1 ut'u_2) \overline{\psi(u_1u_2)} du_1 du_2.\]
If we suppose that $f_0(ut') \neq 0$, then there exist $u_1, u_2 \in U_1$ such that $f(u_1ut'u_2) \neq 0$.  Writing $u_1ut'u_2 = u_1u(t'u_2t'^{-1})t'$, we see that $t' \in \Omega$ and $u_1u(t'u_2t'^{-1}) \in U_1$. This implies that $u \in U_1$, so the change of variables $u_1 \mapsto u_1 u^{-1}$ shows that $f_0(ut')$ is a scalar multiple of 
\[ \psi(u) \int\limits_{U_1 \times U_2} f(u_1 t'u_2) \overline{\psi(u_1u_2)} du_1 du_2. \]
Now we write
\begin{equation*}
    \begin{split}
        \int\limits_{U_1 \times U_2} f(u_1 t'u_2) \overline{\psi(u_1u_2)} du_1 du_2 & = \int\limits_{U_1 \times U_2} f(u_1 (t'u_2t'^{-1})t') \overline{\psi(u_1u_2)} du_1 du_2 \\
        & = \int\limits_{U_1 \times U_2} f(u_1 t') \overline{\psi(u_1(t'u_2^{-1}t'^{-1})u_2)} du_1 du_2 \\
        & = \int\limits_{U_2} \psi(t'u_2t'^{-1}u_2^{-1}) du_2 \int\limits_{U_1} f(u_1t') \overline{\psi(u_1)} du_1.
    \end{split}
\end{equation*}
Identifying $U_M$ with $k$, and $t'$ with $y \in k^{\ast}$, the first integral is
\[ \int\limits_{U_2} \psi(x(y^{-1}-1))dx = 0.\]
This shows that $f_0(ut') = 0$ and completes the proof.  \end{proof}

\subsection{Proof of Theorem 2.2.1}

We now can apply the results of the previous sections to prove the stability result.  Let $\pi_1$ and $\pi_2$ be two supercuspidal representations of $\GL_2(k)$ with the same central character $\omega_{\pi}$.  Let $f_1$ and $f_2$ be matrix coefficients of $\pi_1$ and $\pi_2$ such that $W^{f_i}(e) = 1$.  

Let $\omega$ be a character of $k^{\ast}$, identified with a character of $\GL_2(k)$ through the determinant.  Assume that the central characters of $\pi_1 \otimes \omega$ and $\pi_2 \otimes \omega$ are both ramified.  Then we may apply Shahidi's local coefficient formula (Proposition 2.2.2) for both $C_{\psi}(s, \pi_1 \otimes \omega)$ and $C_{\psi}(s,\pi_2 \otimes \omega)$: there exists an integer $\kappa_{\omega}$, depending on $\omega$, such that 
\begin{equation*}
    \begin{split}
        C_{\psi}(s,\pi_i \otimes \omega)^{-1} = \gamma(40s, \omega_{\pi}^2 \omega^2,\psi)^{-1} \int\limits_R \eta(\det(m)) J_{\varphi_{\kappa}}(n,f_i)\\ (\omega_{\pi} \omega^2)^{-2}( \frac{\frac{1}{2}x_{21} - x_{31}}{x_{21}^2+x_{32}})  |\det(m)|^{10s+\frac{5}{2}} dx_{21}dx_{31}dx_{32}.
    \end{split}
\end{equation*}
for all $\kappa \geq \kappa_{\omega}$.  We have used the fact that if $f_i(m)$ is a matrix coefficient of $\pi_i$, $\omega(\det(m))f_i(m)$ is one of $\pi_i \otimes \omega$.

Now we fix an auxiliary function $f_0 \in \mathscr C_c^{\infty}(M ; \omega_{\pi})$ and apply the results of (2.3).  Since $W^{f_i}(e) = W^{f_0}(e) = 1$ for $i = 1, 2$, we have $W^{f_i - f_0}(e) = 0$.  Proposition 2.2.3 tells us that there exists an $f_i^0 \in \mathscr C_c^{\infty}(M ;\omega_{\pi})$, supported inside the big cell of $M$, such that
\[ J_{\varphi_{\kappa}}(n,f_i) = J_{\varphi_{\kappa}}(n,f_0) + J_{\varphi_{\kappa}}(n,f_i^0) \]
for all $\kappa$ greater than or equal to some $\kappa_0$ depending on $f_1,f_2,$ and $f_0$.  We compute the difference $C_{\psi}(s,\pi_1 \otimes \omega)^{-1} - C_{\psi}(s,\pi_1\otimes \omega)^{-1}$, the common term $J_{\varphi_{\kappa}}(n,f_0)$ cancels.  We obtain the difference of local coefficients as 
\begin{equation*}
    \begin{split}
    \gamma(40s, \omega_{\pi}^2 \omega^2,\psi)^{-1} \int\limits_R \omega(\det(m)) (J_{\varphi_{\kappa}}(n,f_1^0) - J_{\varphi_{\kappa}}(n,f_2^0)) \\ (\omega_{\pi} \omega^2)^{-2}( \frac{\frac{1}{2}x_{21} - x_{31}}{x_{21}^2+x_{32}})  |\det(m)|^{10s+\frac{5}{2}} dx_{21}dx_{31}dx_{32}.
    \end{split}
\end{equation*}
whenever $\kappa \geq \operatorname{Max}\{\kappa_0, \kappa_{\omega}\}$.  By Proposition 2.1.1, there exists a compact open subgroup $H$ of $k^{\ast}$, depending on $f_1^0, f_2^0,$ and $\omega_{\pi}$, such that if $t \in H$, $n = \exp(1,0,x_{21},x_{31},x_{32})$ and $n_1 = \exp(1,0,tx_{21}, tx_{31}, t^2x_{32})$, then $J_{\varphi_{\kappa}}(n,f_i^0) = J_{\varphi_{\kappa}}(n_1,f_i^0)$.  Take $\omega$ to be sufficiently highly ramified so that $\omega_{\pi}^2 \omega^2$ is nontrivial on $H$, and choose a $t \in H$ with $\omega_{\pi} \omega(t^2) \neq 1$.  

In the above integral, make the change of variables $(x_{21}, x_{31}, x_{32}) \mapsto (tx_{21}, tx_{31}, t^2 x_{32})$.  Then $\det(m)$ will change to $t^{-2}\det(m)$, and 
\[ (\omega_{\pi} \omega^2)^{-2}(\frac{\frac{1}{2}x_{21} - x_{31}}{x_{21}^2+x_{32}}) \]
will change to 
\[  \omega_{\pi}(t^2) \omega(t^4)(\omega_{\pi} \omega^2)^{-2}(\frac{\frac{1}{2}x_{21} - x_{31}}{x_{21}^2+x_{32}}). \]
What this shows is that
\[ C_{\psi}(s,\pi_1 \otimes \omega)^{-1} - C_{\psi}(s,\pi_2 \otimes \omega)^{-1} = \omega_{\pi} \omega(t^2) (C_{\psi}(s,\pi_1 \otimes \omega)^{-1} - C_{\psi}(s,\pi_2 \otimes \omega)^{-1} ) \]
and therefore $C_{\psi}(s,\pi_1 \otimes \omega)^{-1}  - C_{\psi}(s,\pi_2 \otimes \omega)^{-1} = 0$.  This completes the proof of Theorem 2.2.1.  

\newtheorem{2_4_1}{Remark}[subsection]

\begin{2_4_1} In the proof of stability for symmetric and exterior square gamma factors given in [CoShTs17], the transfer from $n$ to $m$ under $\dot w_0^{-1}n = mn'\bar{n}$ is quite straightforward.  In our case, and in what we expect for other exceptional Lie groups, the relationship between $n$ and $m$ is more subtle.

Just as in [CoShTs17], the proof of Theorem 2.2.1 used a change of variables in the orbit space integral to conclude that the difference of the local coefficients was zero.  Even though the transfer of $n$ to $m$ is less straightforward in our case, what allowed the proof to work was the fact that the coefficients of $m$ and $\bar{n}$ were quotients of homogenous polynomials in the variables $x_{21}, x_{31}, x_{32}^2$.  If this is the case for other exceptional groups, for example $E_6$, it is possible that the method of proof given here can carry over to that case.

\end{2_4_1}

\end{document}